\documentclass[reqno, 10pt]{amsart}

\usepackage{amsfonts}
\usepackage{graphicx}

\newtheorem{theorem}{Theorem}[section]
\newtheorem{proposition}[theorem]{Proposition}
\newtheorem{lemma}[theorem]{Lemma}

\newtheorem{remark}{Remark}[section]

\newcommand{\no}{\noindent}
\newcommand{\nn}{\nonumber}
\newcommand{\oo}{\infty}

\newcommand{\ra}{\rightarrow}

\makeatletter
\@addtoreset{equation}{section}
\makeatother

\renewcommand\thefigure{\thesection.\@arabic\c@figure}
\renewcommand\thetable{\thesection.\@arabic\c@table}

\bibdata{prob}
\bibstyle{alpha}

\title[Copula for multivariate Extreme Value distributions]
{On the Copula for multivariate Extreme Value distributions}

\author{Marco Aur\'elio Sanfins and Glauco Valle}

\date{}

\address{
\newline
UFRJ - Departamento de m\'etodos estat\'{\i}sticos do Instituto de Matem\'atica.
\newline  Caixa Postal 68530, 21945-970, Rio de Janeiro, Brasil
\newline
e-mail:  \rm \texttt{glauco.valle@im.ufrj.br}
\newline
e-mail:  \rm \texttt{sanfins@dme.ufrj.br}
}

\thanks{}

\begin{document}

\maketitle

\begin{abstract}
We show that all multivariate Extreme Value distributions, which are the possible weak limits of the $K$ largest order statistics of iid samples, have the same copula, the so called K-extremal copula. This copula is described through exact expressions for its density and distribution functions. We also study measures of dependence, we obtain a weak convergence result and we propose a simulation algorithm for the K-extremal copula.
\end{abstract}

\bigskip

\section{Introduction}
\label{sec:intro}

In the study of extremes of iid sequences a question of interest is whether or not the dependence relation among the marginals of the limit distribution of the $K$ largest order statistics relies on the parent distribution function of the sequence. One way to evaluate nonlinear dependence between random variables is through the copula associated to them, this is already discussed in several books as the ones by Joe \cite{joe}, Nelsen \cite{nelsen} and Drouet-Mari and Kotz \cite{marikotz}. In the present paper, we show that every multivariate extreme value distribution, which are the possible weak limits of the $K$ largest order statistics of iid samples, have the same copula called the $K$-extremal copula. From the Extremal Types Theorem, see below, extremal distributions are obtained from linear transformations of one of three basic distributions. We prove that the copula for the three basic types is the $K$-extremal copula, thus all K-dimensional multivariate extremal distribution have the same nonlinear dependence among its marginals. This is not remarkable since the copula for any group of order statistics of an iid sample of size n with continous parent distribution do not denpend on this distribution, see Lemma 6 in \cite{averousgenestkochar}. However, a proper caracterization of the K-extremal copula is relevant as well as their consequences. Our result generalizes the case $K=2$ which was considered in \cite{mendessanfins}.

The $K$-extremal copula is described by its distribution and density functions through exact expressions. We show that the copula of the $K$ largest order statistics of iid sequences with continuous parent distribution converges in distribution to the $K$-extremal copula. We also study the assymptotic behavior of Spearman's rho and Kendall's tau for the first and the $K$ largest order statistics. As a last result, we propose a simulation algorithm to sample from the $K$-extremal copula.

In section \ref{sec:stat} we will present and discuss the results in this paper postponing all the proofs to section \ref{sec:proofs}.

\section{Statements}
\label{sec:stat}

\medskip
Fix an interger $K \ge 2$. For every $n\ge K$, let $M_{1,n}$, ... , $M_{K,n}$ be the $K$ largest order statistics of an iid sample of size $n$ with parent distribution not depending on $n$. The Extremal Types Theorem, see sections 2.2 and 2.3 in \cite{leadbetterlindgrenrootzen} and section 4.2 in \cite{embrechtskluppelbergmikosch}, states that if for some sequences of real numbers $(a_n)_{n=1}^\oo$ and $(b_n)_{n=1}^\oo$ the random variables $a_n M_{1,n} + b_n$ converge in distribution then the random vectors
\begin{equation}
\label{maxnorm}
(a_n M_{1,n} + b_n , ... , a_n M_{K,n} + b_n )
\end{equation}
also converge in distribution. The limit belongs to a family of distributions parametrized by $-\oo < \mu < \oo$, $\sigma > 0$ and $-\oo < \xi < \oo$. For a choice $(\mu,\sigma,\xi)$ of the parameters, the marginals of a limit distribution have
distribution function and density functions given respectively by
\begin{equation}
\label{eq:G}
G_m (z) = \left\{
\begin{array}{ll}
\exp\{- \Lambda(z)\} \sum_{j=0}^{m-1} \frac{\Lambda(z)^j}{j!} \!\!\!\!\! &, \textrm{ if } \xi \left( \frac{z-\mu}{\sigma}\right) > -1 \textrm{ for } \xi \neq 0 \textrm{ or } z \in \mathbb{R} \textrm{ for } \xi = 0  \\
0 &, \textrm{ if } z < \mu - \frac{\sigma}{\xi} \textrm{ for } \xi > 0 \\
1 &, \textrm{ if } z > \mu - \frac{\sigma}{\xi} \textrm{ for } \xi < 0.
\end{array} \right.
\end{equation}
and
\begin{equation}
\label{eq:g}
g_m (z) = \left\{
\begin{array}{ll} \exp\{- \Lambda(z)\} \frac{ \Lambda^\prime(z) \Lambda(z)^{m-1}}{(m-1)!}  \!\!\!\!\! &, \textrm{ if } \xi \left( \frac{z-\mu}{\sigma}\right) > -1 \textrm{ for } \xi \neq 0 \textrm{ or } z \in \mathbb{R} \textrm{ for } \xi = 0  \\
0 &, \textrm{ otherwise,}
\end{array} \right.
\end{equation}
where
$$
\Lambda (z) = \Lambda_{\xi,\mu,\sigma} (z) = \left\{
\begin{array}{ll}
\left[ 1 + \xi \left( \frac{z-\mu}{\sigma}\right) \right]^{-\frac{1}{\xi}} & , \textrm{ if } \xi \neq 0 \\
\exp \left( - \frac{z-\mu}{\sigma} \right) & , \textrm{ if } \xi = 0,
\end{array} \right. \, .
$$
for $m \le 1$. A distribution with distribution function as above is called a Generalized Extreme Value (GEV) distribution which are classified in types I, II and III according respectively to $\xi = 0$, $\xi > 0$ and $\xi < 0$. Note that the function $\Lambda$ is strictly decreasing positive function and satisfies
\begin{eqnarray}
\label{lambda}
\lim_{z \ra -\oo} \Lambda(z) = +\oo \quad \textrm{and} \quad \lim_{z \ra \oo} \Lambda(z) = 0 , \textrm{ if } \xi = 0 \nn \\
\lim_{z \downarrow (\mu - \frac{\sigma}{\xi})} \Lambda(z) = +\oo \quad \textrm{and} \quad \lim_{z \ra \oo} \Lambda(z) = 0 , \textrm{ if } \xi > 0  \\
\lim_{z \ra -\oo} \Lambda(z) = +\oo \quad \textrm{and} \quad \lim_{z \uparrow (\mu - \frac{\sigma}{\xi})} \Lambda(z) = 0 , \textrm{ if } \xi < 0. \nn
\end{eqnarray}

Also by the Extremal Types Theorem, the joint density function $\tilde{g}_K$ of a limiting extreme value distribution for normalized sums of the $K$ largest order statistics of an iid sequence, as in (\ref{maxnorm}), is given by
\begin{equation}
\label{eq:gtil}
\tilde{g}_K (z_1, ... , z_K) = \left\{ \begin{array}{ll}
(-1)^K \exp\{- \Lambda(z_K)\} \, \prod_{j=1}^K \Lambda^\prime (z_j) &, \textrm{ if } (z_1, ..., z_K) \in \Omega_{\xi} \\
0 &, \textrm{ otherwise.}
\end{array} \right. \, ,
\end{equation}
where 
$$
\Omega_{\xi} = \left\{
\begin{array}{ll}
\mathbb{R}^K &\!\!\!\! , \ \textrm{ if } \xi = 0 \\
\{ (z_1, ..., z_K) \in \mathbb{R}^K : z_1> ... > z_K > \mu - \frac{\sigma}{\xi} \} &\!\!\!\! , \ \textrm{ if } \xi > 0 \\
\{ (z_1, ..., z_K) \in \mathbb{R}^K : \mu - \frac{\sigma}{\xi} > z_1> ... > z_K \} &\!\!\!\! , \ \textrm{ if } \xi < 0 .
\end{array}
\right.
$$
A distribution with density as in (\ref{eq:gtil}) for parameters $-\oo < \mu < \oo$, $\sigma > 0$ and $-\oo < \xi < \oo$ is called a Multivariate Generalized Extreme Value (MGEV) distribution.

\smallskip

\begin{remark}
A broader class of stationary sequences of random variables have a MGEV distribution as the assymptotic distribution of the largest maxima. These sequences should satisfy some weak dependence condition. The results can be found for instance in \cite{embrechtskluppelbergmikosch}. 
\end{remark}

\smallskip
Our first result gives an explicity expression for the distribution function associated to the density $\tilde{g}_K$.

\medskip

\begin{proposition}
\label{prop:Gtil}
The distribution function $\tilde{G}_K$ of a limiting extreme value distribution for a normalized vector of the $K$ largest order statistics of iid continuous random variables has the following representation
$$
\tilde{G}_K (z_1, ... , z_K) = H_K (z_1, \min(z_1,z_2), \min(z_1,z_2,z_3), ... , \min(z_1, ... , z_K)) \, ,
$$
for every $(z_1,...,z_K) \in \mathbb{R}^K$, where
$$
H_K (z_1, ... , z_K) = \exp\{- \Lambda(z_K)\} \, J_K ( \Lambda(z_1), ... , \Lambda(z_K)) 
$$
for $\min(z_1,...,z_K) > \mu - \frac{\sigma}{\xi}$, if $\xi > 0$, or for $\min(z_1,...,z_K) < \mu - \frac{\sigma}{\xi}$, if $\xi < 0$, or $(z_1,...,z_K) \in \mathbb{R}^K$, if $\xi=0$, otherwise  $H_K (z_1, ... , z_K) = 0$. The function $J_K : \mathbb{R}_+^K \ra \mathbb{R}_+$ is a polynomial in $K$ variables which is defined by induction by putting $J_1 \equiv 1$ and
$$
J_m (x_1, ... , x_m) =  \sum_{j=0}^{m-1} \frac{x_m^j}{j!} - \sum_{j=1}^{m-1} \frac{x_j^j}{j!} J_{m-j} (x_{j+1}, ... , x_m), \quad \textrm{ for } m\ge 1.
$$
\end{proposition}

\medskip
We can now compute the density of the copula associated to the density $\tilde{g}_K$ of a MGEV distribution, which we call the K-extremal copula and tuns out to not depend on the parameters $\xi$, $\mu$ and $\sigma$.

\medskip

\begin{proposition}
\label{prop:cK}
The density of the copula of a MGEV distribution is given by
\begin{eqnarray}
\label{eq:cK-1} \lefteqn{
c_K (u_1, ... , u_K) =  \left( \prod_{j=1}^{K-1} \frac{d\log\psi_j}{du_j}(u_j) \right) \frac{d\psi_K}{du_K}(u_K)} \\
\label{eq:cK-2} & & = \left( \prod_{j=1}^{K-1} (-1)^{j-1} \psi_j (u_j) \frac{(\log\psi_j(u_j))^{j-1}}{(j-1)!} \right)^{\!\!\! -1} \left( \frac{(-\log\psi_K(u_K))^{K-1}}{(K-1)!} \right)^{\!\!\! -1} \!\!\!\! ,
\end{eqnarray}
for $(u_1,...,u_K) \in (0,1)^K$ such that $u_1 > \psi_2(u_2) > ... > \psi_K(u_K)$, where $\psi_m:(0,1) \ra (0,1)$ is the increasing function that satisfies the following implicit equation
\begin{equation}
\label{eq:psi}
u = \psi_m (u) \sum_{j=0}^{m-1} (-1)^j \frac{(\log\psi_m(u))^j}{j!} \, ,
\end{equation} 
otherwise $c_K (u_1, ... , u_K) = 0$.
\end{proposition}

\smallskip

\begin{remark} 
\label{remark:psi}
The function $\psi_m$ which appears in the expression for the density of the K-extremal copula can be obtained from a MGEV distribution function as $\psi_m (u) = \exp \{ - \Lambda(G_m^{-1}(u) \}$ for every $u\in (0,1)$ and $m \ge 1$.
\end{remark}

\smallskip
Also with the distribution function of the MGEV distribution, it is straightforward to write the distribution function of the K-extremal copula which we present in the next result. 

\medskip

\begin{proposition}
\label{prop:CK}
The copula of a MGEV is given by
$$
C_K (u_1,...,u_K) = \mathcal{H}_K (u_1, r_1(u_1,u_2), r_2(u_1,u_2,u_3), ... , r_{K-1}(u_1, ... , u_K)) \, .
$$
for every $(u_1,...,u_K) \in [0,1]^K$, where
$$
r_{m-1}(u_1,...,u_{m})= \psi_m^{-1}(\psi_l(u_l)) = \psi_l(u_l) \sum_{j=0}^{m-1} (-1)^j \frac{(\log\psi_l(u_l))^j}{j!},
$$
if $\psi_l(u_l) = \min(\psi_1(u_1),...,\psi_m(u_m))$ and for every $(u_1,...,u_K)$ such that $u_1 = \psi_1(u_1) \ge \psi_2(u_2) \ge ... \ge \psi_K(u_K)$
\begin{eqnarray*}
\lefteqn{
\mathcal{H}_K (u_1, ... , u_K) =  \psi_K(u_K) J_K \left( -\log u_1 , -\log \psi_2(u_2), ... , -\log \psi_K (u_K) \right), } \\
& = & u_K - \psi_K(u_K) \sum_{j=1}^{K-1} \frac{(-\log\psi_j(u_j))^j}{j!} J_{K-j}(-\log\psi_{j+1}(u_{j+1}),...,-\log\psi_K(u_K))
\end{eqnarray*}
with $J_m$ defined in the statement of Proposition \ref{prop:Gtil}.
\end{proposition}

\bigskip

By a simple generalization of Lemma 6 in \cite{averousgenestkochar}, we have that the multivariate copula of the $K$ largest order statistics of an iid sample of size n do not depend on the continuous parent distribution of the sample. This copula will be denoted by $\tilde{C}^{(n)}_K$, where $n$ denotes the size of the sample. The next proposition is a convergence result for copulas that has the consequence that for continuous distributions the non-linear dependence structure of the K-largest order statistics of large iid samples is approximatedly captured by the K-extremal copula.

\bigskip

\begin{proposition}
\label{prop:convcop}
The copula $\tilde{C}^{(n)}_K$ converges in distribution to $C_K$ as $n \ra \oo$.
\end{proposition}

\bigskip

From the $K$-extremal copula we can obtain the copula between the l largest and the m largest limiting order statistics for every choice of $l$ and $m$, or between any two marginals of a MGEV distribution. Then we can use these bivariate copulas to obtain measures of dependence as the Spearman's rho and Kendall's tau. For a copula $C$, the Spearman's rho is defined by
$$
12 \int_0^1 \int_0^1 C(u,v) dudv -3 = 12 \int_0^1 \int_0^1 u v dC(u,v) -3
$$
and Kendall's tau by
$$
4 \int_0^1 \int_0^1 C(u,v) dC(u,v) -1 \, .
$$
We are going to study here the behavior of Spearman's rho and Kendall's tau for the first and the $K$th marginals of the $K$-extremal copula in the limit as $K \ra \oo$. We denote these measures respectively by $\rho_K$ and $\tau_K$, $K \ge 2$. Using the convergence result in proposition 2.4, this caracterizes the behavior of these measures for the first and the Kth largest order statistics of large samples with continuous parent distribution. We point out that $\rho_2 = 2/3$ and $\tau_2 = 1/2$ have been obtained \cite{mendessanfins}. For more on measures of dependence of order statistics see \cite{averousgenestkochar} and \cite{chen}. We have the following result.

\begin{proposition}
\label{prop:rhotau}
Both sequences $(\rho_K)$ and $(\tau_K)$ converges to zero as $K \ra \oo$.
\end{proposition}

\bigskip

We now describe a simulation algorithm to generate samples from the K-extremal copula. The method is based on a technique of conditional sampling to sample from multivariate copulas, see for instance Cherubini, Luciano and Vecchiato's book \cite{cherubinilucianovecchiato}. We can resume the procedure with the following steps:
\begin{enumerate}
\item[(i)] Put $C_i(u_1 , u_2 ,..., u_m) = C(u_1 , u_2 ,..., u_m,1,...,1)$ for $m=2,...,K$;
\item[(ii)] Sample $u_1$ from the uniform distribution in $(0,1)$;
\item[(iii)] Sample $u_m$ from the conditional distribution $C_m ( \cdot | u_1,...,u_{m-1})$ for $m=1,...,K$;
\end{enumerate}

We now are going to focus on how to sample $u_k$ from the conditional distribution $C_k ( \cdot | u_1,...,u_{k-1})$. To sample $u_m$ from $C_m (. |  u_1,...,u_{m-1})$, we sample $q$ from $U(0,1)$ and we put $u_m=C^{-1}_m(q|u_1,...,u_{m-1})$. Therefore we should know explicitly $C_m (\cdot | u_1,...,u_{m-1})$. We compute it in the following lemma:

\begin{lemma} \label{conddist} The condicional distribution function of $U_m |(U_1, U_2, ...,U_{m-1})$ when $(U_1, ...,U_{K})$ has distribution given by the K-extremal copula is given by
\begin{equation}
\label{eq:6.9}
C_m (u_m |  u_1,...,u_{m-1}) = \frac{\psi_m ( u_m )}{\psi_{m-1}(u_{m-1})}.
\end{equation}
\end{lemma}

\medskip

If we now put $ q= C_m (u_m |  u_1,...,u_{m-1})$, we have that:
$$
u_m = C^{-1}_m(q|u_1,...,u_{m-1}) = \psi^{-1}_m (q .\psi_{m-1}(u_{m-1}) ) \, .
$$
From definition \ref{eq:psi} we get
$$
u_m =  \psi_m(q .\psi_{m-1}(u_{m-1})) \sum_{j=0}^{m-1} (-1)^j \frac{(\log \psi_m(q .\psi_{m-1}(u_{m-1})))^j}{j!}.
$$
Therefore, we solve numerically $\psi_{m-1}(u_{m-1})$ and then $\psi_m(q .\psi_{m-1}(u_{m-1}))$ to obtain $u_m$.

We plot below a sample of size 200 from the 4-extremal copula.
\begin{figure}[h]
    \centering
        \includegraphics[height=5cm,width=8cm]{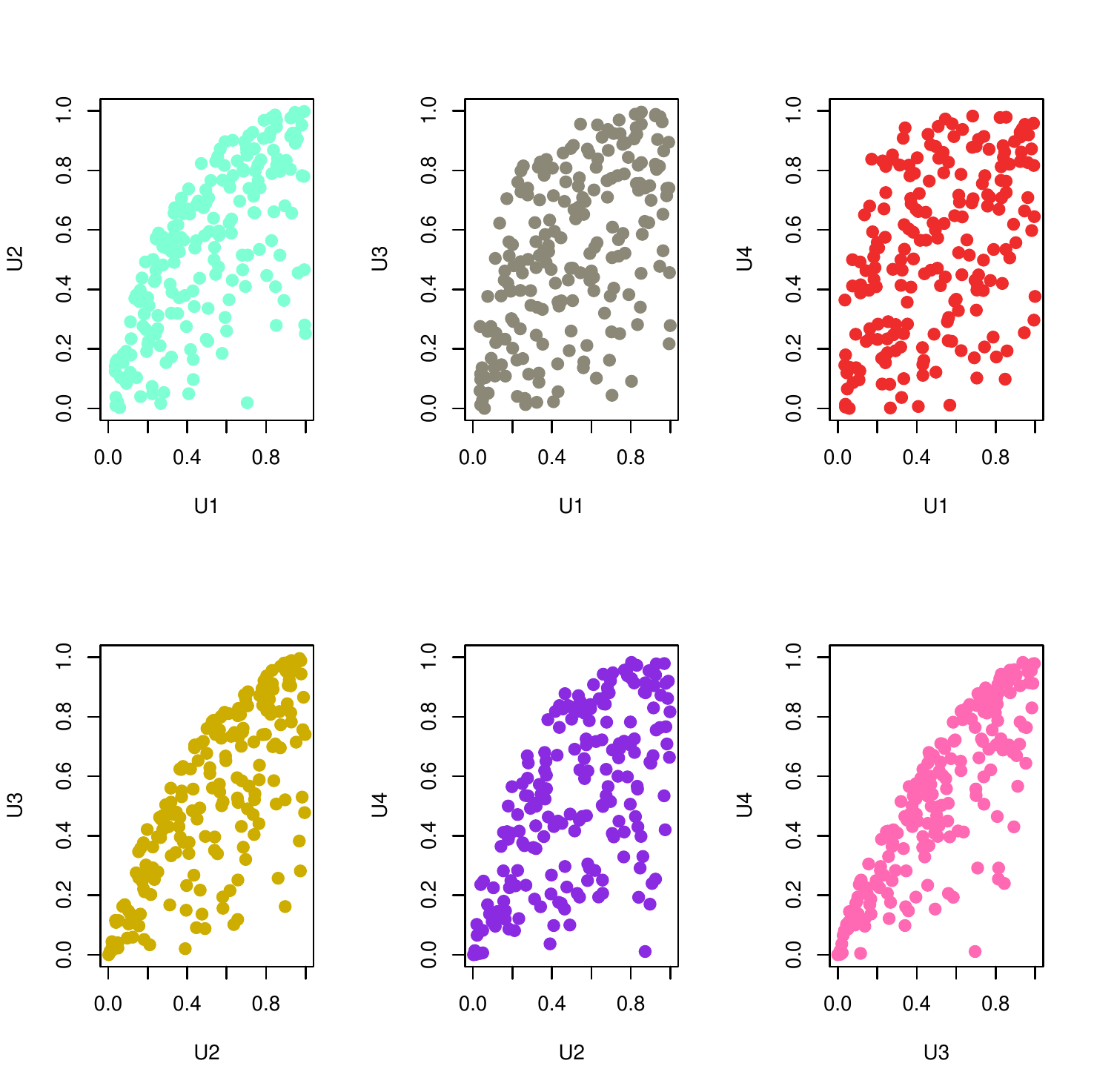}
    \label{fig20}
\end{figure}

\bigskip
\section{Proofs}
\label{sec:proofs}

\medskip

\no  \textbf{Proof of Proposition \ref{prop:Gtil}:}
We show that $\tilde{G}_K$ is a $K$-dimensional distribution function with density given by $\tilde{g}_K$. By the definition of $\tilde{g}_K$, the multiple integral
$$
\int_{-\oo}^{z_1} ... \int_{-\oo}^{z_K} \tilde{g}_K (y_1, ... y_K) \ dy_1...dy_K
$$
is equal to
$$
\int_{-\oo}^{z_1} \int_{-\oo}^{\min(z_1,z_2)} ... \int_{-\oo}^{\min(z_1,...,z_K)} \tilde{g}_K (y_1, ... y_K) \ dy_1...dy_K.
$$
Therefore $\tilde{G}_K (z_1, ... ,z_K) = \tilde{G}_K (z_1, \min(z_1,z_2), ... ,\min(z_1,...,z_K))$. From now on, we suppose that $z_1 > z_2 > ... > z_K$. Then, from (\ref{eq:gtil}),
$$
\tilde{G}_K (z_1, ... ,z_K) = (-1)^K \int_{A_{\xi}}^{z_K} \int_{y_{K}}^{z_{K-1}} ... \int_{y_3}^{z_2} \int_{y_2}^{z_1} \exp\{- \Lambda(y_K)\} \, \prod_{j=1}^K \Lambda^\prime (y_j) \ dy_1...dy_K \, ,
$$
where $A_\xi = \mu - \frac{\sigma}{\xi}$, if $\xi > 0$, and $A_\xi = -\oo$ otherwise. 
Considering the following change of variables in the last integral, $x_j = \Lambda(y_j)$, for $1\le j \le K$, we get the following integral
$$
I_K(w_1,...,w_K) := (-1)^K \int_{w_K}^{+\oo} \int^{x_{K}}_{w_{K-1}} ... \int^{x_3}_{w_2} \int^{x_2}_{w_1} e^{- x_K}  dx_1...dx_K \, ,
$$
where $w_j= \Lambda(z_j)$. To complete the proof, We show by induction that
$$
I_K(w_1,...,w_K) = e^{-w_K} J_K (w_1 , ... , w_K ) \, .
$$
For $K=1$, a simple verification shows that the result holds. Now suppose that it holds for $1 \le K \le L-1$. For $K = L$, we perform the first iterated integral in the expression for $I_K(w_1,...,w_K)$ to obtain that it is equal to
$$
(-1)^K \int_{w_K}^{+\oo} \int^{x_{K}}_{w_{K-1}} ... \int^{x_3}_{w_2} x_2 e^{- x_K}  dx_2...dx_K - w_1 I_{K-1}(w_2,...,w_K) \, .
$$
Then perform the first iterated integral in the first term of the previous expression to obtain
\begin{eqnarray*}
\lefteqn{
(-1)^K \int_{w_K}^{+\oo} \int^{x_{K}}_{w_{K-1}} ... \int^{x_4}_{w_3} \frac{x_3}{2} e^{- x_K}  dx_3...dx_K - } \\
& & \qquad \qquad \qquad \qquad - \frac{w_2}{2} I_{K-2}(w_3,...,w_K) - w_1 I_{K-1}(w_1,...,w_K) \, .
\end{eqnarray*}
Following recursively this procedure we get
$$
I_K(w_1,...,w_K) = e^{-w_K} \sum_{j=0}^{m-1} \frac{w_K^j}{j!} - \sum_{j=1}^{m-1} \frac{w_j^j}{j!} I_{K-j}(w_{j+1},...,w_K).
$$
By the definition of $J_K$ and the induction hypotheses we complete the proof. $\square$

\bigskip

\no \textbf{Proof of Proposition \ref{prop:cK}:}
Let us fix a limiting extreme value distribution function $\tilde{G}_K$. We have that
$$
c_K(u_1,...,u_K) = \frac{\tilde{g}_K(G_1^{-1}(u_1),...,G_K^{-1}(u_K))}{\prod_{j=1}^K g_j(G_j^{-1}(u_j))} \, .
$$
Therefore we just apply formulas (\ref{eq:g}) and (\ref{eq:gtil}) to obtain that $c_K(u_1,...,u_K)$ is equal to
$$
\left( \prod_{j=1}^{K-1} \exp\{-\Lambda(G_j^{-1}(u_j))\} \frac{\Lambda(G_j^{-1}(u_j))^{j-1}}{(j-1)!} \right)^{\!\!\! -1} \left( \frac{\Lambda(G_K^{-1}(u_K))^{K-1}}{(K-1)!} \right)^{\!\!\! -1} .
$$
From this formula, if we put $\psi_m (u) = \exp \{ - \Lambda(G_m^{-1}(u) \}$ we get (\ref{eq:cK-2}) in the statement. Now (\ref{eq:psi}) is a direct consequence of the explicit formulas for the distribution function $G_m$ given in (\ref{eq:G}). 

It remains to verify (\ref{eq:cK-1}). If we derive both sides of (\ref{eq:psi}), we get that
\begin{eqnarray*}
1 & = & \left( \sum_{j=0}^{m-1} (-1)^j \frac{(\log\psi_m)^{j}}{(j)!} - \sum_{j=0}^{m-2} (-1)^j \frac{(\log\psi_m)^{j}}{(j)!} \right) \frac{d\psi_m}{du} \\
& = & (-1)^{m-1} \frac{(\log\psi_m)^{m-1}}{(m-1)!} \frac{d\psi_m}{du},
\end{eqnarray*}
which implies that
\begin{equation}
\label{derpsi}
\frac{d\psi_m}{du} = (-1)^{m-1} \left( \frac{(\log\psi_m)^{m-1}}{(m-1)!} \right)^{\!\!\! -1}
\end{equation}
and
\begin{equation}
\label{derlogpsi}
\frac{d\log\psi_m}{du} = (-1)^{m-1} \left( \psi_m \frac{(\log\psi_m)^{m-1}}{(m-1)!} \right)^{\!\!\! -1} \!\! .
\end{equation}
From (\ref{derpsi}), (\ref{derlogpsi}) and (\ref{eq:cK-2}) we arrive at (\ref{eq:cK-1}).
$\square$

\bigskip

\no \textbf{Proof of Proposition \ref{prop:CK}:} Let us fix a limiting extreme value distribution function $\tilde{G}_K$. Then the distribution function of the K-extremal copula is given by
$$
C_K(u_1,...,u_K) = \tilde{G}_K (G_1^{-1}(u_1),...,G_K^{-1}(u_K)) 
$$
for every $(u_1,...,u_K) \in [0,1]^K$ which by Proposition \ref{prop:Gtil} is equal to
$$
H_K (G^{-1}_1(u_1), \min(G^{-1}_1(u_1),G^{-1}_2(u_2)) , ... , \min(G^{-1}_1(u_1),...,G^{-1}_K(u_K))).
$$
By the definition of $H_K$, monotonicity and the expression for $\psi_m$ in remark \ref{remark:psi}, see also the proof of Proposition \ref{prop:cK}, the previous expression is equal to
$$
\min_{1\le l \le K} ( \psi_l(u_l) ) \ J_K \left( -\log u_1 , -\log \min_{l=1,2} ( \psi_l(u_l) ), ... , -\log \min_{1\le l \le K} (\psi_l(u_l) ) \right).
$$
Using the definition of $r_m$ in the statement, write the above expression as
$$
\psi_K (r_K(u_1,...,u_m)) \ \ J_K \left( -\log u_1 , -\log \psi_2(r_2(u_1,u_2)), ... , -\log \psi_K (r_K(u_1,...,u_m)) \right),
$$
which completes the proof. $\square$

\bigskip

\no \textbf{Proof of Proposition \ref{prop:convcop}:} Let $M_{1,n}$, ... , $M_{K,n}$ be the K-largest order statistics of a sample of size $n$ with a given continuous parent distribution function $F$ which belongs to the domain of atraction of a GEV distribution. This means that there exists $(a_n)_{n=1}^{+\oo}$ and $(b_n)_{n=1}^{+\oo}$ sequences of real numbers such that the random vector
$$
(a_n M_{1,n} + b_n , ... , a_n M_{K,n} + b_n )
$$
converges in distribution to some $\tilde{G}_K$ which is MGEV distribution. By invariance concerning composition with affine transformations the copula associated to $(M_{1,n} , ... ,M_{K,n} )$ and $(a_n M_{1,n} + b_n , ... , a_n M_{K,n} + b_n )$ is $\tilde{C}^{(n)}_K$ independently of $F$.

Let $F_{j,n}$ be the distribution function of $a_n M_{j,n} + b_n$. Therefore, if we define the function $V_n(x_1,...,x_K) =(F_{1,n}(x_1) , ... , F_{K,n}(x_K) )$, $(x_1,...,x_K) \in \mathbb{R}^n$ then
\begin{equation}
\label{maxnormcopula}
V_n(a_n M_{1,n} + b_n,...,a_n M_{K,n} + b_n)
\end{equation}
has the distribution of the copula $\tilde{C}^{(n)}_K$.

The K-extremal copula has the distribution of $V(Y_1,...,Y_K)$, where $V(x_1,...,x_K) =(G_{1}(x_1) , ... , G_{K}(x_K) )$, $(x_1,...,x_K) \in \mathbb{R}^n$. By Theorem 5.1 in \cite{billingsley}, (\ref{maxnormcopula}) converges in distribution to the K-extremal copula if $V_n$ converges uniformly to $V$ on compact intervals, but this is a consequence of P\'olyas's Theorem which implies that $F_{j,n}$ converges uniformly to $G_j$ since the last is absolutely continuous. $\square$

\bigskip

\no \textbf{Proof of Proposition \ref{prop:rhotau}:} We shall prove through estimates on exact expressions that $\rho_K \ra 0$. The analogous result can be applied to $\tau_K$ since $\rho_K \ge \tau_K \ge 0$. This last assertion can be verified through Theorem 5.1 of Fredricks and Nelsen in \cite{fredricksnelsen}. Indeed, according to their terminology, for two order statistics, the largest is always left-tail decreasing and smallest is right-tail increasing.

Applying directly the definition we can write $(\rho_K +3)/12$ as
\begin{equation}
\label{eq:intrhoA}
\int_0^1 \int_{\psi^{-1}_{K-1}(\psi_K(u_K))}^1 ... \int_{\psi^{-1}_{2}(\psi_3(u_3))}^1 \int_{\psi_2(u_2)}^1 u_1 \, u_K \, c_K(u_1, ... ,u_K) \, du_1...du_K.
\end{equation}
which we are going to show that converges to $1/4$ as $K \ra \oo$ resulting in $\rho_K \ra 0$. By (\ref{eq:cK-1}) the previous iterated integral can be rewritten as
$$
\int_0^1 \int_{\psi^{-1}_{K-1}(\psi_K(u_K))}^1 \!\!\!\!\!\! ... \int_{\psi^{-1}_{2}(\psi_3(u_3))}^1 \int_{\psi_2(u_2)}^1 \!\!\!\!\!\! u_1 \,  u_K  \left( \prod_{j=1}^{K-1} \frac{d\log\psi_j}{du_j}(u_j) \right) \! \frac{d\psi_K}{du_K}(u_K) \, du_1...du_K.
$$
By induction in $1\le m \le K-1$, we show that
$$
\int_{\psi^{-1}_{m}(\psi_{m+1}(u_{m+1}))}^1 ... \int_{\psi^{-1}_{2}(\psi_3(u_3))}^1 \int_{\psi_2(u_2)}^1 u_1 \, \prod_{j=1}^{m} \frac{d\log\psi_j}{du_j}(u_j) \, du_1...du_m.
$$
is equal to
\begin{equation}
\label{eq:intrhoB}
(-1)^m \left[ \psi_{m+1}(u_{m+1}) - \sum_{j=0}^{m-1} \frac{(\log \psi_{m+1}(u_{m+1}))^j}{j!} \right] \, .
\end{equation}
Indeed, $\psi_1$ is the identity function in $(0,1)$ and therefore
$$
\int^1_{\psi_2(u_2)} u_1 \frac{d\log \psi_1}{du_1}(u_1) du_1 = (-1)[\psi_2(u_2) -1]\, .
$$
Now suppose that (\ref{eq:intrhoB}) holds for some $1 \le l \le K-2$ then
$$
(-1)^l \left[ \psi_{l+1}(u_{l+1}) - \sum_{j=0}^{l-1} \frac{(\log \psi_{l+1}(u_{l+1}))^j}{j!} \right] \frac{d\log \psi_{l+1}}{du_{l+1}}(u_{l+1}) \, .
$$
is equal to
$$
(-1)^l \frac{d}{du_{l+1}} \! \left( \psi_{l+1}(u_{l+1}) - \sum_{j=1}^{l} \frac{(\log \psi_{l+1}(u_{l+1}))^j}{j!} \right)
$$
and, since $\psi_{l+1}(1)=1$, integrating on $u_{l+1}$ over the interval $(\psi^{-1}_{l+1} (\psi_{l+2}(u_{l+2})),1)$ we obtain that (\ref{eq:intrhoB}) holds for $m= l+1$.

Therefore the integral in (\ref{eq:intrhoA}) is equal to 
$$
\int_0^1 u \, \frac{d\psi_K}{du}(u) (-1)^{K-1} \left[ \psi_{K}(u) - \sum_{j=0}^{K-2} \frac{(\log \psi_{K}(u))^j}{j!}  \right] du \, .
$$
Put $v=\psi_K(u)$, $u \in (0,1)$ and uses the power series expansion
$$
v = \sum_{j=0}^\oo \frac{\log(v)^j}{j!}
$$ 
to write the previous integral as
$$
(-1)^{K-1} \int_0^1 \psi^{-1}_K (v) \left( \sum_{j=K-1}^\oo \frac{\log(v)^j}{j!} \right) dv \, .
$$
Another change of variables and (\ref{eq:psi}) allows us to write the integral in (\ref{eq:intrhoA}) as 
$$
(-1)^{K-1} \sum_{l=0}^{K-1} \sum_{j= K-1}^\oo \frac{(-1)^j}{j!l!} \int_0^{+\oo} y^{l+j} e^{-2y} dy 
$$
which, since
$$
\int_0^{+\oo} y^{l+j} e^{-2y} dy = \frac{(l+j)!}{2^{l+j+1}} \, ,
$$
can be rewritten as 
$$
 (-1)^{K-1} \sum_{l=0}^{K-1} \sum_{j= K-1}^\oo (-1)^j \Big( \!\! \begin{array}{c} l+j \\ l \end{array} \!\! \Big) \frac{1}{2^{l+j+1}} \, .
$$
We finish the proof showing that
$$
\lim_{K \ra \oo} \, \left\{ (-1)^{K-1} \sum_{l=0}^{K-1} \sum_{j= K-1}^\oo  (-1)^j  \Big( \!\! \begin{array}{c} l+j \\ l \end{array} \!\! \Big) \frac{1}{2^{l+j}} \right\} = \frac{1}{2} \, .
$$
From this point we suppose that $K$ is odd, for $K$ even the proof is similar with few sign changes. The left hand side term in the previous convergence statement is equal to
\begin{equation}
\label{eq:intrhoC}
\sum_{l=0}^{K-1} \sum_{j= K-1}^\oo  \Big( \!\! \begin{array}{c} l+j \\ l \end{array} \!\! \Big) \frac{1}{2^{l+j}} \, - \sum_{l=0}^{K-1} \sum_{j= \frac{K-1}{2}}^\oo  \Big( \!\! \begin{array}{c} l+2j+1 \\ l \end{array} \!\! \Big) \frac{1}{2^{l+2j}} \, ,
\end{equation}
Now apply the identities
$$
\Big( \!\! \begin{array}{c} l+2j \\ 0 \end{array} \!\! \Big) = 1 \quad \textrm{and} \quad  
\Big( \!\! \begin{array}{c} l+2j+1 \\ l \end{array} \!\! \Big) = \Big( \!\! \begin{array}{c} l+2j \\ l -1 \end{array} \!\! \Big) + \Big( \!\! \begin{array}{c} l+2j \\ l  \end{array} \!\! \Big) \, , \ \textrm{for } l\ge 1,
$$
to write the second term in (\ref{eq:intrhoC}) as
$$
\sum_{j=K-1}^\oo  \Big( \!\! \begin{array}{c} 2j+1 \\ K-1 \end{array} \!\! \Big) \frac{1}{2^{2j+1}} - \sum_{l=0}^{K-1} \sum_{j= K-1 }^\oo  \Big( \!\! \begin{array}{c} l+j \\ l \end{array} \!\! \Big) \frac{1}{2^{l+j}} \, .
$$
Therefore (\ref{eq:intrhoC}) is equal to
$$
\sum_{j=K-1}^\oo  \Big( \!\! \begin{array}{c} 2j+1 \\ K-1 \end{array} \!\! \Big) \frac{1}{2^{2j+1}} 
$$
which is 
$$
\sum_{j=2K}^\oo  \Big( \!\! \begin{array}{c} j-1 \\ K-1 \end{array} \!\! \Big) \frac{1}{2^{j}} + \sum_{j=K-1}^\oo \Big( \!\! \begin{array}{c} 2j+1 \\ K-1 \end{array} \!\! \Big) \left( 1 - \frac{2j+2}{2(2j-K+3)} \right) \frac{1}{2^{2j+2}}.
$$
Let $Y$ be a random variable with negative binomial distribution with parameters $K$ and $1/2$. Then the second term in the sum above is equal to
$$
\mathbb{E} \left[ \left( 1 - \frac{Y}{2(Y-K+1)} \right) \textrm{I} \{Y \textrm{ even}, \ Y \ge 2K \} \right] \, ,
$$
which is bounded above by
\begin{eqnarray*}
\lefteqn{\mathbb{E} \left[ \left( 1 - \frac{Y}{2(Y-K+1)} \right) \textrm{I} \{ 2K \le Y \ge 2K + K^{\frac{3}{4}} \} \right] 
+ \mathbb{P} (  Y \ge 2K + K^{\frac{3}{4}} )} \\
& & \le \left( 1 - \frac{2 + K^{-\frac{1}{4}}}{2 + 2K^{-\frac{1}{4}}+2K^{-1})} \right) + 
\mathbb{P} \left(  \frac{Y - \mathbb{E}[Y]}{\sqrt{2K}} \ge \frac{K^{\frac{1}{4}}}{\sqrt{2}} \right)
\end{eqnarray*}
that goes to zero as $K \ra \oo$ by the central limit theorem.

Therefore the limit of (\ref{eq:intrhoC}) as $K \ra \oo$ is the same as the limit of
$$
\sum_{j=2K}^\oo  \Big( \!\! \begin{array}{c} j-1 \\ K-1 \end{array} \!\! \Big) \frac{1}{2^{j}} 
$$
which is the probability that a negative binomial distribution with parameters $K$ and $1/2$ takes a value greater or equal to $2K$. This probability converges to $1/2$ again by the Central Limit Theorem. $\square$

\bigskip

\no \textbf{Proof of Lemma \ref{conddist}:} Let $(U_1, U_2, ...,U_K)$ be a random vector whose distribution function is $C$. The conditional distribution of $U_m$ given $U_1, U_2, ...,U_{m-1}$ has distribution function
\begin{eqnarray}
\label{eq:6.8}
C_m (u_m | u_1,...,u_{m-1}) \!\!\!\!\!  & = & \!\!\!\!\! \mathbb{P} (U_m \le u_m | U_1 = u_1 ,..., U_{m-1} = u_{m-1}) \nonumber \\[9pt]
                             & = & \!\!\!\!\! \frac{ \left( \frac{\partial^{m-1} C_m(u_1,...,u_m)}{ \partial u_1,...,\partial u_{m-1}} \right) }{ \left( \frac{\partial^{m-1} C_{m-1}(u_1,...,u_{m-1})}{\partial u_1,...,\partial u_{m-1}} \right) }
\end{eqnarray}
for every $m=2,...,k$.

We first deal with the numerator in (\ref{eq:6.8}) which by the formula in Proposition \ref{prop:CK} can be written as
\begin{eqnarray}
\frac{\partial^{m-1} \left[- \psi_m (u_m) \sum_{j=1}^{m-1} \frac{-log(\psi_j (u_j))^j}{j!} J_{m-j} (-\log \psi_{j+1} (u_{j+1}) , ...,-\log \psi_m (u_m )) \right]}{\partial u_1 ... \partial u_{m-1}} \nn \, .
\end{eqnarray}
If we remove the terms that do not depend on all the variables $u_1 , ... , u_{m-1}$, we obtain that
the last partial derivative is equal to
\begin{eqnarray}
\label{eq:6.9}
\frac{\partial^{m-1} \left[ - \psi_m (u_m ) \prod_{j=1}^{m-1} (-log (\psi(u_j ) ) )\right]}{\partial u_1 ... \partial u_{m-1}} \, .
\end{eqnarray}
Using that
$$
\frac{d \log \psi_m}{du} = (-1)^{m-1} \left( \psi_m \frac{(\log \psi_m)^{m-1}}{{(m-1)}!}\right)^{-1} \, ,
$$
we obtain that (\ref{eq:6.9}) is equal to
\begin{eqnarray}
\label{eq:6.11}
(-1)^{m}\psi_m (u_m) (-1)^{m-1} \prod_{j=1}^{m-1} (-1)^{j-1} \left( \psi_j (u_j) \frac{log(\psi_j (u_j))^{j-1}}{(j-1)!} \right)^{-1} \, .
\end{eqnarray}
Now we consider the denominator in (\ref{eq:6.8}) which is equal to the density function of the (m-1)-extremal copula. Hence it is equal to
\begin{equation}
\label{eq:6.12}
\left( \prod_{j=1}^{m-2} (-1)^{j-1}  \psi_j (u_{j})\frac{(\log \psi_j (u_j))^{j-1})}{(j-1)!}  \right)^{-1}\left(  -\frac{(\log \psi_{m-1} (u_{m-1}))^{m-2}}{(m-2)!}\right)^{-1} \, .
\end{equation}

Finally replace the expressions in (\ref{eq:6.11}) and (\ref{eq:6.12}) respectively in the numerator and denominator  in (\ref{eq:6.8}) to obtain that
$$
C_m (u_m |  u_1,...,u_{m-1}) = \frac{\psi_m ( u_m )}{\psi_{m-1}(u_{m-1})}. \quad \square
$$

\bigskip



\begin{thebibliography}{99}

\bibitem{averousgenestkochar} J. Averous, C. Genest, S. Kochar: On dependence structure of order statistics, Journal of Mult. Anal. 94, 159-171 (2005)

\bibitem{billingsley} P. Billingsley: \emph{Convergence of Probability Measures}, John Wiley $\&$ Sons (1968)

\bibitem{cherubinilucianovecchiato} U. Cherubini, E. Luciano, W. Vecchiato: \emph{Copula Methods in Finance}. Wiley Finance: Chichester. (2004)

\bibitem{fredricksnelsen} G. Fredricks, R. Nelsen: On the relation between Spearman's rho and Kendall's tau for pairs of continuous random variables, Journal of Stat. Plan. and Inf. 137, 2143-2150 (2007)

\bibitem{marikotz} D. Drouet-Mari, S. Kotz: \emph{Correlation and dependence}, Imperial College Press, London (2001).

\bibitem{embrechtskluppelbergmikosch} P. Embrechts, C. Kluppelberg, T. Mikosch: \emph{Modelling Extremal Events for Insurance and Finance}, Applications of Mathematics 33, Springer

\bibitem{joe} H. Joe: \emph{Multivariate models and dependence concepts}, Chapman $\&$ Hall (1997)

\bibitem{leadbetterlindgrenrootzen} M. R. Leadbetter; G. Lindgreen; H. Rootzen: \emph{Extremes and related properties of random sequences and processes}, Springer-Verlag New York (1983)

\bibitem{mendessanfins} B. Mendes, M. Sanfins: The limiting copula of the two largest order statistics of iid samples, \emph{Brazilian Jounal of prob. and Statisc.} (2007)

\bibitem{nelsen} R. Nelsen: {\it An introduction to copulas}, Springer series in statistics, $2^o$ edition (2007).

\bibitem{chen} Y. Chen: A note on the relantionship between Spearman's $\rho$ and Kendall's $\tau$ for extreme order statistics, \emph{Journal of Stat. Plan. and Inf.} 137, 2165-2171 (2007)

\end{thebibliography}
\end{document}